\documentclass[a4paper,11pt]{article}
\usepackage[dvips]{graphicx}
\usepackage{here}
\usepackage[top=25truemm,bottom=25truemm,left=25truemm,right=25truemm]{geometry}
\usepackage{amsthm}
\usepackage{amsmath,amssymb}
\usepackage{color}
\usepackage{cite}
\usepackage{authblk}

\begin{document}

\title{Exploitation by asymmetry of information reference in coevolutionary learning in prisoner's dilemma game}

\author{Yuma Fujimoto$^a$, Kunihiko Kaneko$^b$}

\affil{Department of Evolutionary Studies of Biosystems, School of Advanced Sciences, SOKENDAI, Japan.$^a$\\
Department of Basic Sciences, The University of Tokyo, Japan.$^b$}
\affil{{\rm fujimoto\_yuma@soken.ac.jp}$^a$}

\maketitle

\begin{abstract}
Mutual relationships, such as cooperation and exploitation, are the basis of human and other biological societies. The foundations of these relationships are rooted in the decision making of individuals, and whether they choose to be selfish or altruistic. How individuals choose their behaviors can be analyzed using a strategy optimization process in the framework of game theory. Previous studies have shown that reference to individuals' previous actions plays an important role in their choice of strategies and establishment of social relationships. A fundamental question remains as to whether an individual with more information can exploit another who has less information when learning the choice of strategies. Here we demonstrate that a player using a memory-one strategy, who can refer to their own previous action and that of their opponent, can be exploited by a reactive player, who only has the information of the other player, based on mutual adaptive learning. This is counterintuitive because the former has more choice in strategies and can potentially obtain a higher payoff. We demonstrated this by formulating the learning process of strategy choices to optimize the payoffs in terms of coupled replicator dynamics and applying it to the prisoner's dilemma game. Further, we show that the player using a memory-one strategy, by referring to their previous experience, can sometimes act more generous toward the opponent's defection, thereby accepting the opponent's exploitation. Mainly, we found that through adaptive learning, a player with limited information usually exploits the player with more information, leading to asymmetric exploitation.
\end{abstract}

\section{Introduction}
Cooperation, defection, and exploitation are important relationships that universally appear in biological and social systems. While cooperating, individuals are altruistic and achieve benefits for the entire group. In defection, they behave selfishly for their own benefit, which results in demerits for all. In exploitation, selfish individuals receive benefit at the expense of altruistic others. The choice of strategy, i.e., selfish or altruistic behavior, is important in establishing social relationships. Individuals, based on their abilities, sophisticate their strategies through their experiences. Generally, people's ability to choose the best strategies differ. These differences in ability can affect how cooperation is established between people. Now the following question arises: Do individuals with higher abilities exploit those with lower abilities or vice versa?

\begin{figure}[tbhp]
\begin{center}
\includegraphics[width=0.8\linewidth]{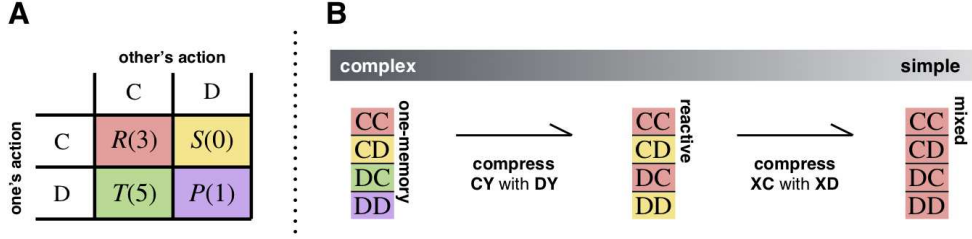}
\caption{(A) Payoff matrix for the single PD game. A single game has four results: CC (red), CD (yellow), DC (green), and DD (purple). (B) memory-one, reactive, and mixed strategy classes. The memory-one class can refer to all four previous results. The reactive class can only refer to an opponent's result, compressing CC and DC (colored in the same red), and CD and DD (yellow). Furthermore, the mixed class compresses all the previous results into one (all colored in red).}
\label{F01}
\end{center}
\end{figure}

Game theory is a mathematical framework for analyzing such individual decision-making of strategies \cite{vonNeumann2007}. Everyone has a strategy for choosing given actions and receives a reward based on their chosen actions. In particular, the prisoner's dilemma (PD) game (see Fig.~\ref{F01}-(A)) has been used extensively to investigate how people act when competing for benefits. Each of two players chooses either cooperation (C) or defection (D), depending on their own strategies. Accordingly, the single game has four results, namely CC, CD, DC, and DD, where the left and right symbols (C or D) indicate the action taken by oneself and that of the opponent, respectively. The benefit of them is given by $R$, $S$, $T$, and $P$. The property of PD demands $T>R>P>S$: Each player receives a larger benefit for choosing D, which may lead to DD, but CC is more beneficial than DD. How to avoid falling into mutual defection, i.e., DD, has been a significant issue.

In this study, we assume the iteration of games, and then each player can refer to the actions of the previous game and change the own choice depending on the observed actions.  We consider that a player can change their next action based on the previous actions of the two players, i.e., CC, CD, DC, and DD at the maximum. Thus, a player with more detailed information about a previous round has a higher ability to choose their own optimal action. This ability to observe actions of their own and their opponents is seen in reality, such as intention recognition \cite{Premack1978,Saxe2007,Lurz2011,Han2011,Han2015,Moniz2015,Fujimoto2019a}. Furthermore, in social systems, from the bacterial community to international war, the reference to past actions might play an important role in establishing interpersonal relationships \cite{Axelrod1984}. A representative example is tit-for-tat (TFT) strategy \cite{Axelrod1981,Axelrod1988}, in which the player observes and mimics the other's previous action. The monumental and the numerous following studies showed that the players with TFT strategies are selected in the optimization process and establish cooperation. This strategy is classified as a reactive strategy \cite{Nowak1990a,Nowak1990b,Nowak1992,Imhof2010,Zhang2011,Baek2016} as the player chooses their actions by referring only to their opponent's previous action. Another type of strategy, the memory-one strategy \cite{Baek2016,Nowak1993,Brauchli1999,Kraines2000,Iliopoulos2010,Stewart2014,Hilbe2017}, is introduced when a player refers to both their own previous action as well as their opponent's. The memory-one strategy includes not only the TFT but the Win-Stay-Lose-Shift (WSLS) strategy \cite{Nowak1993}, which generates cooperation even under the error in choice of action. Consequently, it is expected that a player who refers to more information will succeed in receiving a larger benefit.

Although such reference of past actions plays a crucial role to establish interpersonal relationships, how the difference in information for it between the players affects in the player's benefit is unclear yet. Evolution of strategy in a single population \cite{Nowak1990a,Nowak1990b,Nowak1992,Imhof2010,Zhang2011,Baek2016,Nowak1993,Brauchli1999,Kraines2000,Iliopoulos2010,Stewart2014,Hilbe2017} does not generate the difference in payoffs among the players in principle, where the strategy is selected within the same group. Evolution in multi-population \cite{Traulsen2006,Wang2013,Jin2014,Xia2018,Liu2019,Takesue2021} (i.e., intra-group selection of strategies by the inter-group game) and learning among multi-agent \cite{Macy1991,Macy1998,Macy2002a,Macy2002b,Sandholm1996,Taiji1999,Masuda2009,Fujimoto2019b} can generate the exploitation, but so far most studies have focused only on the emergence of cooperation. Recently, the exploitation has been studied as to ``symmetry'' breaking of the payoffs of players \cite{Jin2014,Fujimoto2019b}, where only the game between the same class of strategies is assumed. Here, we revised the coupled replicator model in the previous studies of multi-population evolution \cite{Sigmund1998} and multi-agent learning \cite{Borgers1997,Sato2002,Sato2003} so that the player can refer to their own previous actions and those of their opponent and update their strategy accordingly within a class of strategies. In particular, we focused on the reactive and memory-one classes of strategies, as they are basic and have been studied extensively. We then investigated whether players using the memory-one strategy would win the game against opponents using the reactive strategy, by utilizing the extra information provided through the observation of their own previous actions.

The remainder of this paper is organized as follows. In \S~\ref{S2}, we formulate the learning dynamics for various strategy classes. In addition, we confirm that in iterated games against an opponent's fixed strategy, the strategy with the higher ability obtains a larger payoff in equilibrium. In \S~\ref{S3}, we introduce an example of mutual learning between memory-one and reactive strategies. Then, we demonstrate that the memory-one class, i.e., the player with the higher ability, is, counterintuitively, one-sidedly exploited by the reactive one. In \S~\ref{S4}, we analyze how this exploitation is achieved and elucidate that the ability to reference one's own actions leads to generosity and leaves room for exploitation. Finally, in \S~\ref{S5}, we show that high-ability players are generally exploited because of their generosity, independent of the strategy class and payoff matrix.

\section{Formulation of learning dynamics of strategies}
\label{S2}
\subsection{Formulation of class and strategy}
\label{S2-0}
Before formulating the learning process, we mathematically define the strategy and class. Recall that a single game can have one of four results: CC, CD, DC, and DD. First, a player using a memory-one strategy can refer to their own and the opponent's previous actions and respond with a different action to each result of the previous game. Thus, the player has four independent stochastic variables $x_1$, $x_2$, $x_3$, and $x_4$ as the probabilities of choosing C regardless of the outcome of the previous game. Thus, the memory-one class is defined as the possible set of such memory-one strategies, which is denoted as $\{x_1,x_2,x_3,x_4\}\in[0,1]^4$. The TFT and Win-Stay-Lose-Shift \cite{Nowak1993,Posch1999,Imhof2007,Amaral2016} strategies are examples of the emergence of cooperation as $x_1=x_3=1$, $x_2=x_4=0$, and $x_1=x_4=1$ and $x_2=x_3=0$.

Second, a player using the reactive strategy can only refer to the opponent's action and therefore cannot distinguish between CC and DC (CD and DD). Thus, the strategy is given by two independent variables $x_{13}$ and $x_{24}$, where the former (latter) is the probability of choosing C where the previous result was either CC or DC (CD or DD). Therefore, the reactive class is defined as $\{x_{13},x_{24}\}\in[0,1]^2$. Here, the notation $x_{13}$ (variable for CC and DC) clearly indicates the integration of $x_1$ (CC) and $x_3$ (DC) from the memory-one class. Indeed, all the strategies in the reactive class are included in the memory-one class, as one can set $x_1=x_3=x_{13}$ and $x_2=x_4=x_{24}$ for all $x_{13}$ and $x_{24}$. Thus, the above TFT strategy can be represented as $x_{13}=1$ and $x_{24}=0$, whereas the Win-Stay-Lose-Shift strategy cannot be represented. This makes it clear that the memory-one class is more complex than the reactive one, because the former can include all the strategies of the latter. The ordering of complexity can be defined as a class that includes all the strategies of the other is more complex.

Third, in the classical mixed strategy \cite{Nash1950}, a player stochastically chooses their action without referencing any actions from the previous game. Thus, the strategy controls only one variable $x_{1234}\in[0,1]$, which is the probability of choosing C. This class is the least complex of the three classes.

\subsection{Analysis of repeated game}
\label{S2-1}
In this section, we analyze a repeated game under the condition that the strategy of both players are fixed. We define ${\bf p}:=(p_{\rm CC},p_{\rm CD},p_{\rm DC},p_{\rm DD})^{\mathrm{T}}$ as the probabilities that (CC, CD, DC, DD) are played in the present round, and the result of the next round is calculated by ${\bf p}'={\bf Mp}$ with
\begin{equation}
	{\bf M}:=\left(\begin{array}{cccc}
		x_{1}y_{1} & x_{2}y_{3} & x_{3}y_{2} & x_{4}y_{4} \\
		x_{1}\bar{y_{1}} & x_{2}\bar{y_{3}} & x_{3}\bar{y_{2}} & x_{4}\bar{y_{4}} \\
		\bar{x_{1}}y_{1} & \bar{x_{2}}y_{3} & \bar{x_{3}}y_{2} & \bar{x_{4}}y_{4} \\
		\bar{x_{1}}\bar{y_{1}} & \bar{x_{2}}\bar{y_{3}} & \bar{x_{3}}\bar{y_{2}} & \bar{x_{4}}\bar{y_{4}} \\
	\end{array}\right).
\label{E5-001}
\end{equation}

When none of the strategy variables, $x_n$ and $y_n$ for $n\in\{1,\cdots,4\}$ are $0$ or $1$, the repeated game has only one equilibrium state ${\bf p}_{\mathrm{e}}:=(p_{\rm CCe},p_{\rm CDe},p_{\rm DCe},p_{\rm DDe})^{\rm T}$. Here, we can directly compute ${\bf p}_{\rm e}$ as
\begin{equation}
\begin{split}
	&p_{\rm CCe}=k\{(x_{4}+(x_{3}-x_{4})y_{2})(y_{4}+(y_{3}-y_{4})x_{2})-x_{3}y_{3}(x_{2}-x_{4})(y_{2}-y_{4})\}, \\
	&p_{\rm CDe}=k\{(x_{4}+(x_{3}-x_{4})y_{4})(\bar{y_{2}}-(y_{1}-y_{2})x_{1})-x_{4}\bar{y_{1}}(x_{1}-x_{3})(y_{2}-y_{4})\}, \\
	&p_{\rm DCe}=k\{(\bar{x_{2}}-(x_{1}-x_{2})y_{1})(y_{4}+(y_{3}-y_{4})x_{4})-\bar{x_{1}}y_{4}(x_{2}-x_{4})(y_{1}-y_{3})\}, \\
	&p_{\rm DDe}=k\{(\bar{x_{2}}-(x_{1}-x_{2})y_{3})(\bar{y_{2}}-(y_{1}-y_{2})x_{3})-\bar{x_{2}}\bar{y_{2}}(x_{1}-x_{3})(y_{1}-y_{3})\}. \\
\end{split}
\label{E5-002}
\end{equation}
Here the coefficient $k$ is determined by the normalization of the probabilities $p_{\rm CCe}+p_{\rm CDe}+p_{\rm DCe}+p_{\rm DDe}=1$.

\subsection{Learning dynamics of memory-one class}
\label{S2-2}
We next consider adaptive learning from past experiences of repeated games. For instance, we assume that the probability of CC being the result of the previous round is $p_{\rm CCe}$. Then, the next action being C (D) would have the probability of $x_{1}$ ($\bar{x_{1}}$). Here, we define $u_{\rm CC(C)}$ ($u_{\rm CC(D)}$) as the benefit that the player gains by performing action C (D). First, the time evolution of $x_1$ is assumed to depend on the amount of experience: the previous game's result and the action in the present one must be CC and C, respectively. Thus, $\dot{x_1}$ is proportional to $p_{\rm CCe}x_{1}$. Second, $\dot{x_1}$ also depends on the benefit of action C and thus is proportional to $u_{\rm CC(C)}-(x_{1}u_{\rm CC(C)}+\bar{x_{1}}u_{\rm CC(D)})$. To summarize, we get
\begin{equation}
	\dot{x_{1}}=x_{1}\bar{x_{1}}p_{\rm CCe}(u_{\rm CC(C)}-u_{\rm CC(D)}).
\label{E5-003}
\end{equation}

Next, we compute $u_{\rm CC(C)}$ and $u_{\rm CC(D)}$. When the previous game's result and the present self-action are CC and C, respectively, the present state is given by ${\bf p}={\bf p}_{\rm CC(C)}:=(y_{1},\bar{y_{1}},0,0)^{\rm T}$. If ${\bf p}_{\rm CC(C)}\neq {\bf p}_{\rm e}$ holds, then the state gradually relaxes to equilibrium with the repetition of the game. Thus, $u_{\rm CC(C)}$ is the total payoff generated by ${\bf p}_{\rm CC(C)}$ until equilibrium is reached, which is given by
\begin{equation}
	u_{\rm CC(C)}=\sum_{t=0}^{\infty}{\bf M}^t({\bf p}_{\rm CC(C)}-{\bf p}_{\rm e})\cdot {\bf u}.
\label{E5-004}
\end{equation}
Here, we define ${\bf u}:=(R,S,T,P)^{\rm T}$ as the vector for the payoff matrix. By contrast, when the previous game's result and the present self-action are CC and D, respectively, the present state is given by ${\bf p}={\bf p}_{\rm CC(D)}:=(0,0,y_{1},\bar{y_{1}})^{\rm T}$. Then, $u_{\rm CC(D)}$ is computed in the same way using
\begin{equation}
	u_{\rm CC(D)}=\sum_{t=0}^{\infty}{\bf M}^t({\bf p}_{\rm CC(D)}-{\bf p}_{\rm e})\cdot {\bf u}.
\label{E5-005}
\end{equation}
By substituting Eqs.~\ref{E5-004} and \ref{E5-005} into Eq.~(\ref{E5-003}), we can write the learning dynamics of $x_{1}$ as
\begin{equation}
	\dot{x_{1}}=x_{1}\bar{x_{1}}p_{\rm CCe}\sum_{t=0}^{\infty}{\bf M}^t({\bf p}_{\rm CC(C)}-{\bf p}_{\rm CC(D)})\cdot {\bf u},
\label{E5-006}
\end{equation}
using only strategy variables, $x_{n}$ and $y_{n}$ for $n\in\{1,\cdots,4\}$, and the payoff variables $(T,R,P,S)$. Similarly, we can derive the time evolution of the other strategy variables $x_{2}$, $x_{3}$, and $x_{4}$.

Eq.~(\ref{E5-006}) appears complicated at first glance, but it can be simplified as
\begin{equation}
	\dot{x_{1}}=x_{1}\bar{x_{1}}\frac{\partial {\bf p}_{\rm e}}{\partial x_{1}}\cdot {\bf u},
\label{E5-007}
\end{equation}
as will be shown in subsection \S~\ref{S2-3}. The same equations hold for the learning of $x_2$, $x_3$, and $x_4$, and for opponent player. Notably, this equation reproduces the original coupled replicator model \cite{Sigmund1998,Borgers1997}.

\subsection{Detailed calculation of learning dynamics}
\label{S2-3}
In this section, we prove that Eqs.~\ref{E5-006} and \ref{E5-007} are equivalent. First, because ${\bf p}_{\rm e}$ is the equilibrium state for the Markov transition matrix ${\bf M}$, we obtain
\begin{equation}
	({\bf E}-{\bf M}){\bf p}_{\rm e}=0.
\label{E5-008}
\end{equation}
Here, by a perturbation of strategy $\delta x_{1}$, the equilibrium state $\delta {\bf p}_{\rm e}$ changes accordingly. By substituting $x_{1}\rightarrow x_{1}+\delta x_{1}$ and ${\bf p}_{\rm e}\rightarrow{\bf p}_{\rm e}+\delta {\bf p}_{\rm e}$ into Eq.~(\ref{E5-001}), we obtain:
\begin{equation}
\begin{split}
	&({\bf E}-{\bf M})\delta {\bf p}_{\rm e}-p_{\rm CCe}({\bf p}_{\rm CC(C)}-{\bf p}_{\rm CC(D)})\delta x_{1}=0\\
	&\Rightarrow \sum_{t=0}^{\infty} {\bf M}^t ({\bf E}-{\bf M}) \frac{\partial {\bf p}_{\rm e}}{\partial x_{1}}=p_{\rm CCe}\sum_{t=0}^{\infty} {\bf M}^t ({\bf p}_{\rm CC(C)}-{\bf p}_{\rm CC(D)}), \\
	&\Rightarrow \left({\bf E}-{\bf M}^{\infty}\right)\frac{\partial {\bf p}_{\rm e}}{\partial x_{1}}=p_{\rm CCe}\sum_{t=0}^{\infty} {\bf M}^t ({\bf p}_{\rm CC(C)}-{\bf p}_{\rm CC(D)}), \\
	&\Rightarrow \frac{\partial {\bf p}_{\rm e}}{\partial x_{1}}=p_{\rm CCe}\lim_{n\rightarrow\infty}\sum_{t=0}^n {\bf M}^t ({\bf p}_{\rm CC(C)}-{\bf p}_{\rm CC(D)}).
\end{split}
\label{E5-009}
\end{equation}
Here, we use ${\bf M}^{\infty} \partial {\bf p}_{\rm e}/\partial x_{1}=0$ because ${\bf M}$ has only one eigenvector for which the eigenvalue is $1$ as the preservation of probability (i.e., $p_{\rm CCe}+p_{\rm CDe}+p_{\rm DCe}+p_{\rm DDe}=1$).

Eq.~(\ref{E5-009}) not only provides a simple representation of the time evolution but is also useful for the numerical simulation of Eq.~(\ref{E5-006}). The right-hand side of Eq.~(\ref{E5-009}) requires an approximate numerical calculation for all eigenvectors of ${\bf M}$. By contrast, the left-hand side demands only the information on the equilibrium state ${\bf p}_{\rm e}$, which is analytically given by Eq.~(\ref{E5-002}).

\subsection{Learning dynamics of other strategies}
\label{S2-4}
In the previous sections, we formulated the learning dynamics of memory-one class strategies against another within the same class. In this section, we consider other cases in which both learned and learning players adopt the reactive class.

First, we consider a case in which a learned player uses the reactive class, and the learning player uses the memory-one class. In this case, the learning is easily given by
\begin{equation}
	\dot{x_{n}}(\{x_{1},\cdots,x_{4}\},\{y_{13},y_{24}\})=\dot{x_{i}}|_{y_{1}=y_{3}=y_{13},y_{2}=y_{4}=y_{24}}
\label{E5-010}
\end{equation}
for $n\in\{1,\cdots,4\}$ because the learned reactive player's strategy is constrained by $y_{1}=y_{3}=y_{13}$ and $y_{2}=y_{4}=y_{24}$.

Second, we consider a case in which the learning player uses the reactive class, and the learned player uses the memory-one class. In this case, the learning player's strategy is given by $(x_{13},x_{24})$. Recall that the learning speed of our model depends on the amount of experience. Because the frequency of observing the opponent's previous action, C, is the total of both CC and DC, the time evolution of $x_{13}$ is the sum of $x_{1}$ and $x_{3}$. Thus, we get
\begin{equation}
\begin{split}
	&\dot{x_{13}}(\{x_{13},x_{24}\},\{y_{1},\cdots,y_{4}\})=(\dot{x_1}+\dot{x_3})|_{x_{1}=x_{3}=x_{13},x_{2}=x_{4}=x_{24}},\\
	&\dot{x_{24}}(\{x_{13},x_{24}\},\{y_{1},\cdots,y_{4}\})=(\dot{x_2}+\dot{x_4})|_{x_{1}=x_{3}=x_{13},x_{2}=x_{4}=x_{24}}.\\
\end{split}
\label{E5-011}
\end{equation}

\section{Numerical result for learning}
\label{S3}
\subsection{One-sided learning against a fixed strategy}
\label{S3-1}
\begin{figure}[tbhp]
\begin{center}
\includegraphics[width=0.4\linewidth]{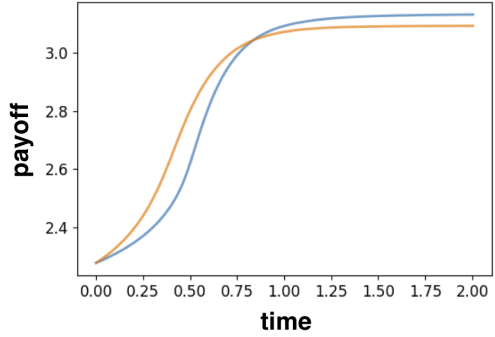}
\caption{Payoff of the memory-one (blue) and reactive (orange) classes over time when learning with an opponent with a fixed strategy, whose strategy is a reactive one; $y_{13}=0.9$ and $y_{24}=0.1$. The payoff is an average of a large number (10000) of initial conditions for $x_i$, which are chosen randomly. The horizontal axis denotes time on a scale of $\log(t+1)$. The rise of the payoff of the reactive class is larger than that of the memory-one class, but finally, the memory-one class has a larger payoff.}
\label{F02}
\end{center}
\end{figure}

Before investigating the game between the memory-one and reactive classes, we first study the learning of the memory-one and reactive classes against the other fixed strategy. Fig.~\ref{F02} shows the time series of each strategy's payoff based on the learning dynamics. The payoffs of both classes monotonically increase their payoffs over time because the opponent's strategy is fixed. However, there are two major differences between the two classes in the way the payoff increases.

First, the reactive class learns faster than the memory-one class. This is because the reactive class is a compressed version of the memory-one model, as the constraints $x_{1}=x_{3}$ and $x_{2}=x_{4}$ are postulated, that is, the learning in the cases of CC and DC (CD and DD) are integrated. Recall that the change in strategy is optimized based on the empirical data sampled through the played games. In the reactive class, the number of strategy variables is fewer; therefore, quick optimization can be achieved, as shown in Eq.~(\ref{E5-011}).

Second, the memory-one strategy gains a larger payoff in equilibrium than the reactive one. This is simply because the memory-one strategy contains a reactive strategy. Accordingly,
\begin{equation}
	\max_{x_{1},x_{2},x_{3},x_{4}}{\bf p}_{\mathrm{e}}(\{x_{1},x_{2},x_{3},x_{4}\},{\bf y})\cdot {\bf u}\ge \max_{x_{13},x_{24}}{\bf p}_{\mathrm{e}}(\{x_{13},x_{24}\},{\bf y})\cdot {\bf u},
\label{E5-012}
\end{equation}
is derived for all the opponent's fixed strategies ${\bf y}$.

\subsection{Mutual learning between memory-one and reactive classes}
\label{S3-2}
In \S~\ref{S3-1}, we considered one-sided learning, where a player dynamically optimizes the strategy against their opponent's fixed strategy. In this section, we consider mutual learning, where both players optimize their strategies as the opponent's strategy continues to change.

We excluded the mixed class in this study because the results of matches with the mixed class are trivial. A player with a mixed class must use the same action independently of the opponent's previous actions. Thus, the opponent always receives a higher pay off by choosing D according to the payoff matrix of the PD. When the opponent's choice is always D, the best choice of the mixed class is also D. Thus, only the pure DD will result in equilibrium. Therefore, the mixed class can establish neither cooperative nor exploitative relationships.

Therefore, we consider only the game between players using the memory-one and reactive classes. In our model, the dynamics of players' strategies are deterministic. Thus, the equilibrium state is uniquely determined by the initial values of strategies ${\bf x}$ and ${\bf y}$. Here, we take sampling over the initial conditions. A match between each pair of classes was evaluated using a sufficiently large number of initial conditions of ${\bf x}$ and ${\bf y}$. In each sample, the initial values of the strategy were assumed to be given randomly. In other words, when the former player takes a memory-one (reactive) class, the strategy is randomly chosen from $\{x_1,x_2,x_3,x_4\}\in[0,1]^4$ ($\{x_{13},x_{24}\}\in[0,1]^2$).

\begin{figure}[tbhp]
\begin{center}
\includegraphics[width=0.8\linewidth]{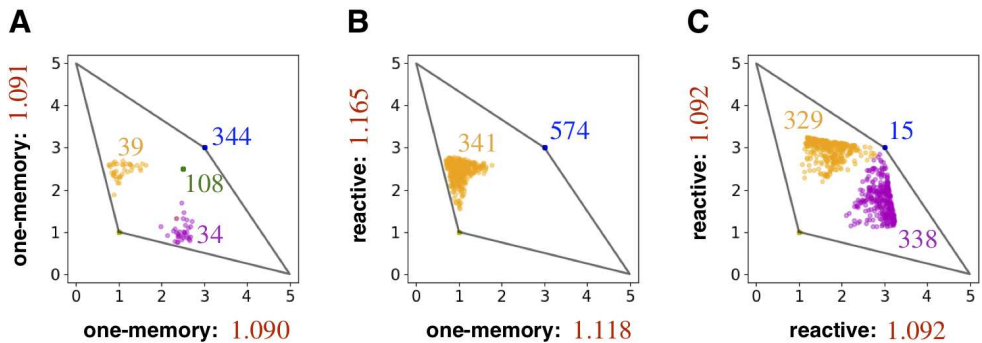}
\caption{Result of matches (A) between two memory-one classes, (B) between a memory-one and reactive class, and (C) between two reactive classes. Each dot indicates one sample, which is given by a random initial state, and $10000$ samples from different initial conditions are shown in each panel. In all panels, the horizontal and vertical lines indicate the player's and the opponent's payoffs. The numbers highlighted in red along the axes represent the average payoff over $10000$ samples. The solid black line indicates the region of the possible set of payoffs. The blue dots indicate CC, where $p_{\rm CCe}=1\Leftrightarrow u_{\rm e}=v_{\rm e}=3$ holds. The yellow dots indicated DD, where $p_{\rm CCe}=0\Leftrightarrow u_{\rm e}=v_{\rm e}=1$ holds. The green dots represent halfway cooperation, where CD and DC are achieved alternately and $p_{\rm CDe}=p_{\rm DCe}=0.5\Leftrightarrow u_{\rm e}=v_{\rm e}=2.5$ holds. The purple (orange) dots indicates the exploitative relationships with $p_{\rm CDe}<p_{\rm DCe}\Leftrightarrow u_{\rm e}>v_{\rm e}$ ($p_{\rm CDe}>p_{\rm DCe}\Leftrightarrow u_{\rm e}<v_{\rm e}$). The numbers in the panels indicate the number of samples that achieved the indicated relationship.}
\label{F03}
\end{center}
\end{figure}

Fig.~\ref{F03} shows the final state of mutual learning for three matches: (A) between two memory-one classes, (B) between the memory-one and reactive classes, and (C) between two reactive classes. (Note that the last case represented in (C) was already studied in [Fujimoto2019].) Here, recall that mutual cooperation satisfies $p_{\rm CCe}=1$, mutual defection satisfies $p_{\rm DDe}=1$, and exploitation satisfies $p_{\rm CDe}\neq p_{\rm DCe}$.

First, we studied the matches between the same classes, represented in (A) and (C). In these matches, exploitation with $p_{\rm CDe}\neq p_{\rm DCe}$ can be in equilibrium. In other words, asymmetry is permanently established between the players depending on their initial strategies, even though both deterministically improve their own strategies to receive a larger payoff. Notably,  this asymmetry emerges symmetrically between the players when using the same class. In this case, the number of samples that satisfied $p_{\rm CDe}>p_{\rm DCe}$ was equal to those that satisfied $p_{\rm CDe}<p_{\rm DCe}$. In (C), i.e., the match between the reactive classes, the equilibrium exists as multiple fixed points with $p_{\rm CDe}\neq p_{\rm DCe}$ (see \cite{Fujimoto2019b} for the detailed analysis). By contrast, each exploitative state in (A) permanently oscillated to form a limit cycle, whereas the temporal averages of $p_{\rm CDe}$ and $p_{\rm DCe}$ are not equal. There is an infinite number of limit cycles, one of which is achieved depending on the initial conditions. A detailed analysis of these limit cycles will be explained in the next section \S~\ref{S4}.

The heterogeneous match (B) between the memory-one and reactive classes has the same exploitative states as match (A). However, the most remarkable difference here is that in this exploitation, the reactive class can receive a larger payoff in the match with the memory-one class, and the reverse never occurs. In other words, only the one-sided exploitation from the reactive class to the memory-one class emerges, regardless of the initial conditions. This result appears paradoxical, when one notes that the memory-one class has more information for the strategy choices and is indeed in a more advantageous position than the reactive one in equilibrium when the other player's strategy is fixed, as already confirmed in \S~\ref{S3-1}. We will discuss the origin of this unintuitive result in Section \S~\ref{S4}.

\section{Emergence of oscillatory exploitation}
\label{S4}
\subsection{Analysis of exploitation}
\label{S4-1}
\begin{figure}[tbhp]
\begin{center}
\includegraphics[width=0.8\linewidth]{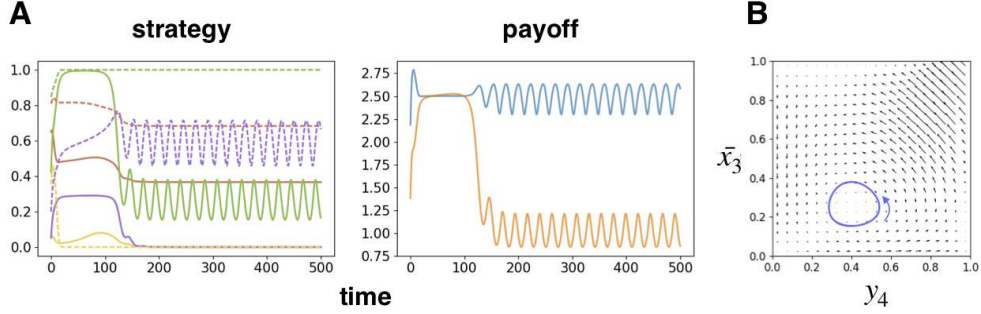}
\caption{The analysis of an exploitative relationship. Panel (A) indicates the trajectory of strategy and payoff for one sample of exploitation by player 1 of player 2. Here, the red, yellow, green, and purple lines indicate $x_1$ ($y_1$), $x_2$ ($y_2$), $x_3$ ($y_3$), and $x_4$ ($y_4$), respectively, whereas the solid and broken lines indicate the exploiting ($x$) and exploited ($y$) player's strategies, respectively. The blue and orange lines indicate the payoffs of the exploiting and exploited player, respectively. Panel (B) shows the vector field for $x_3$ and $y_4$ in the exploitative state. The blue circle is the trajectory of the sample in (A).}
\label{F04}
\end{center}
\end{figure}

We first analyzed the exploitation between the memory-one classes, but the analysis is also applicable to the case between memory-one and reactive classes. An example of the trajectory of strategies $x_i$ and $y_i$ during exploitation is shown in Fig.~\ref{F04}. For all cases, the exploiting player's strategy satisfies
\begin{equation}
\begin{split}
	\begin{array}{ll}
		x_{1}\mathrm{: fixed\ point}, &x_{2}\rightarrow 0,\\
		x_{3}\mathrm{: oscillates}, &x_{4}\rightarrow 0.\\
	\end{array}
\end{split}
\label{E5-015}
\end{equation}
On the other hand, the exploited opponent's strategy satisfies
\begin{equation}
\begin{split}
	\begin{array}{ll}
		y_{1}\mathrm{: fixed\ point}, & y_{2}\rightarrow 0,\\
		y_{3}\rightarrow 1, & y_{4}\mathrm{: oscillates}.\\
	\end{array}
\end{split}
\label{E5-016}
\end{equation}
Here, note that $x_1$ and $y_1$ are neutrally stable, and the asymptotic value continuously varies with the initial condition. Assuming that $x_{2}=x_{4}=y_{2}=0$ and $y_{3}=1$ and inserting them into Eq.~(\ref{E5-002}), the possibility vector ${\bf p}$ satisfies
\begin{equation}
	{\bf p}\propto (0, y_{4}x_{3}, y_{4}, \bar{x_{3}})^{\rm T}.
\label{E5-017}
\end{equation}
This equation leads to $p_{\rm DCe}>p_{\rm CDe}$, which proves that player 1 always receives a larger payoff than player 2. By inserting this into the learning dynamics in Eq.~(\ref{E5-007}), we obtain
\begin{equation}
\begin{split}
	&\dot{x_{1}}=\dot{y_{1}}=0,\\
	&\dot{x_{3}}=x_{3}((T-2P+S)-(T-S)y_{4})\frac{y_{4}\bar{x_{3}}}{(1+y_{4}-x_{3}+y_{4}x_{3})^2},\\
	&\dot{y_{4}}=\bar{x_{4}}((T-P)x_{3}-(P-S))\frac{y_{4}\bar{x_{3}}}{(1+y_{4}-x_{3}+y_{4}x_{3})^2}.
\end{split}
\label{E5-018}
\end{equation}
These equations indicate that $x_{1}$ and $y_{1}$ are neutral, as expected. This is because the players do not experience CC in this oscillatory equilibrium of exploitation and do not have the chance to change $x_{1}$ or $y_{1}$ by learning. Now, the two-variable dynamics $x_{3}$ and $y_{4}$ are obtained, which leads to oscillation.

The oscillatory dynamics for $(x_3,y_4)$ follow the Lotka-Volterra type equation. Eq.~(\ref{E5-018}) has an infinite number of periodic solutions (cycles), and the cycle that is reached depends on the initial strategy. An example of a trajectory is presented in Fig.~\ref{F04}-(B). In Eq.~(\ref{E5-018}), there is only one fixed point, $x_{3}^*$ and $y_{4}^*$, which is given by
\begin{equation}
	(x_{3}^*,y_{4}^*)=\left(\frac{P-S}{T-P}, \frac{T-2P+S}{T-S}\right).
\label{E5-019}
\end{equation}
Then, the players' expected payoffs, $u_{\rm e}^*$ and $v_{\rm e}^*$, are given by
\begin{equation}
	(u_{\rm e}^*,v_{\rm e}^*)=\left(\frac{T+S}{2},P\right).
\label{E5-020}
\end{equation}
This linear stability analysis shows that this fixed point is neutrally stable, that is, the fixed point is not a focus but a center, as in the original Lotka-Volterra equation. Indeed, the time evolution of $(x_3(t),y_4(t))$ has a conserved quantity, given by
\begin{equation}
	H(x_{3},y_{4})=-(P-S)(\log x_{3}+2\log \bar{y_{4}})+(T-P)x_{3}-(T-S)y_{4},
\label{E5-021}
\end{equation}
which is determined by the initial condition and preserved.

Furthermore, this Lotka-Volterra type oscillation provides an explanation of the exploitation we observed. The original Lotka-Volterra equation shows the prey-predator relationship, where the predator increases its own population by sacrificing the prey population. Here $\bar{x_3}$ ($y_4$) represents the exploiter's defection (cooperation on the exploited side), which is a selfish (altruistic) action in the PD. Fig.~\ref{F04}-(B) shows that $\bar{x_3}$ is larger when $x_4$ is larger. In other words, the exploiting side learns to use the selfish action with the altruistic action of the exploited. This result means that the exploiting one increases its own payoff at the expense of the exploited side. Thus, the oscillation of the exploitative relationship is interpreted as a prey-predator relationship.

\subsection{Mechanism of one-sided exploitation: self-reference leads to generosity}
\label{S4-2}
In the previous sections, we mathematically showed how the exploitative relationship is maintained. Next, we intuitively interpret the strategies in Eqs.~\ref{E5-015} and \ref{E5-016}, which implies that exploitation emerges between the narrow-minded and generous players. Here, we also focus on why the memory-one class is exploited one-sidedly by the reactive one.

Before analyzing Eqs.~\ref{E5-015} and \ref{E5-016}, we  present the well-known tit-for-tat (TFT) and two related strategies in Table~\ref{T5-01}. In the TFT strategy, the player deterministically responds with C to the opponent's previous C, and with D to the opponent's D. A more generous strategy \cite{Nowak1992} is that the player accepts the opponent's D and probabilistically responds with C. In contrast, in a more narrow-minded strategy \cite{Fujimoto2019b}, the player betrays the opponent's C probabilistically and responds with D. In contrast to the TFT strategy that was adopted to represent the emergence of symmetric cooperation, the generous and narrow-minded TFT strategies represent asymmetric exploitation. However, these strategies do not refer to previous self-actions. Here, we analyze the exploiting and exploited strategies with the generous and narrow-minded TFT, under the constraint that the previous actions of the self are C and D.
\begin{table}[tbhp]
\begin{center}
\begin{tabular}{ccc}
	\hline
	Strategy & Action to opponent's C & Action to opponent's D \\
	\hline
	TFT & Deterministic C & Deterministic D \\
	Generous TFT & Deterministic C & Probabilistic C \\
	Narrow-minded TFT & Probabilistic D & Deterministic D \\
	\hline
\end{tabular}
\caption{Summary of TFT and related strategies. The TFT strategy deterministically responds with C (D) to the opponent's C (D). Thus, if the player does not refer to their own action , the strategy is represented by a reactive class with $x_{13}=1$ and $x_{24}=0$. The generous TFT's response to the opponent's C is the same as the original TFT but probabilistically cooperates with the opponent's D. This strategy is represented by $x_{13}=1$ and $x_{24}>0$. Finally, a narrow-minded TFT probabilistically defects the opponent's C, which is different from the original TFT. This strategy is represented by $x_{13}<1$ and $x_{24}=0$.}
\label{T5-01}
\end{center}
\end{table}

As seen in Eq.~(\ref{E5-015}), the exploiting player uses the narrow-minded TFT strategy in both cases where their previous action was C and D. In other words, the strategy is characterized by $x_1>0$ and $x_2=0$ for a self C, and $x_3>0$ and $x_4=0$. Therefore, CC never occurs in exploitative equilibrium, where $x_1>0$ can be arbitrarily chosen. Thus, a player has the potential to use the exploiting strategy without referring to their own action, as $x_{13}=x_1=x_3>0$ and $x_{24}=x_2=x_4=0$. Similarly, the exploited player also uses the narrow-minded TFT for a previous self  C, that is, $y_1>0$ and $y_2=0$. However, for a self D, the player uses the generous TFT, that is, $y_3=1$ and $y_4>0$. Thus, the exploited player refers to their own action: the player cannot take $x_{13}=x_1=x_3$ and $x_{24}=x_2=x_4$. Although this additional reference to self-action enriches the player's choice of strategy, the player instead tends to be more generous to the opponent's defection and accepts the opponent's exploitation.

In the above, we consider only the equilibrium states given by Eqs.~\ref{E5-015} and \ref{E5-016}. However, the question remains whether a memory-one class acquires generosity and accepts exploitation during the transient learning process. We attempt to answer this question in the following three steps: First, we classify three equilibrium states in the game between reactive classes: mutual defection, cooperation, and exploitation. Second, by assuming that one of the players adopt a memory-one class instead of the reactive one under these three equilibrium states, we discuss whether the player changes the strategy for each of the three equilibrium states. Third, we consider the one-sided learning process by the above memory-one class under the equilibrium states and examine if the memory-one side's learning increases the opponent's payoff by acquiring generosity.

\begin{figure}[tbhp]
\begin{center}
\includegraphics[width=0.6\linewidth]{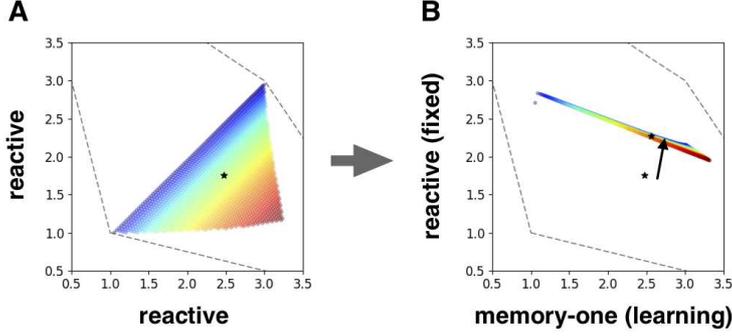}
\caption{In the left panel, payoffs of all the equilibrium states between reactive classes are plotted excluding symmetry. The color map indicates the degree of exploitation. The black star indicates the average value. The right panel shows the payoffs when the exploiting player changes the class into a memory-one class and learns the strategy one-sidedly. The black stars indicate the change in the average payoff. Thus, the learning of the memory-one side leads to the benefit of the reactive opponent rather than for the memory-one itself.}
\label{F05}
\end{center}
\end{figure}

{\bf Step 1.} As shown in Fig.~\ref{F03}-(C), there are three cases of equilibria in a match between reactive classes; mutual defection (yellow dots), mutual cooperation (blue dots), exploitation (purple and orange dots). Here, all possible equilibria for exploitation are shown in Fig.~\ref{F05}.

{\bf Step 2.} Next, we assume that in each of the above states, one player adopts a memory-one class instead of a reactive one. Here, the player can refer to their own actions, after which the state can be unstable. First, for cases of mutual defection and cooperation, the equilibrium states are stable even if one player alternatively adopts a memory-one class. This is easily explained by the analytical result that all the equilibrium states between reactive classes are completely included in those among the memory-one and reactive classes. By contrast, in case of exploitation, the memory-one class receives a larger payoff by releasing constraints $x_1=x_3=x_{13}$ and $x_2=x_4=x_{24}$.

{\bf Step 3.} Finally, we consider one-sided learning by a memory-one class under the exploitation case. Since we assume that the strategy of the reactive opponent is fixed, the memory-one side receives a larger payoff. For all states, the memory-one class releases the constraints $x_1=x_3$ and $x_2=x_4$ and learns to be $x_1=x_2=0$ and $x_3=x_4=1$, which represents an asymptotic relationship to the exploited strategy of Eq.~(\ref{E5-016}). Then, the learning of the memory-one side increases the opponent's benefit much more than the increase in one's own benefit (see Fig.~\ref{F05}). This result shows that the memory-one class becomes generous toward the opponent's defection by also referring to their own actions.

Thus, we have shown that memory-one class becomes generous for the opponent's defection in equilibrium. If the learning feedback from the opponent reactive class is considered, the one-sided exploitation is generated as seen in Fig.~\ref{F03}-(B).

\section{Generality of the result}
\label{S5}
\subsection{Generality over different classes}
\label{S5-1}
We have demonstrated that the reference to own actions lead to generosity toward the opponent's defection when comparing memory-one and reactive classes. Recall that the reactive class identifies both CC and DC, and CD and DD as the same. Therefore, we can consider two intermediate classes between the memory-one and reactive classes: the player only compresses either CC with DC or CD with DD. These strategies refer to the opponent's action completely but refer to the self-action only when the opponent's action is C or D. In this section, we study the learning dynamics of such extended classes.

\begin{figure}[tbhp]
\begin{center}
\includegraphics[width=0.7\linewidth]{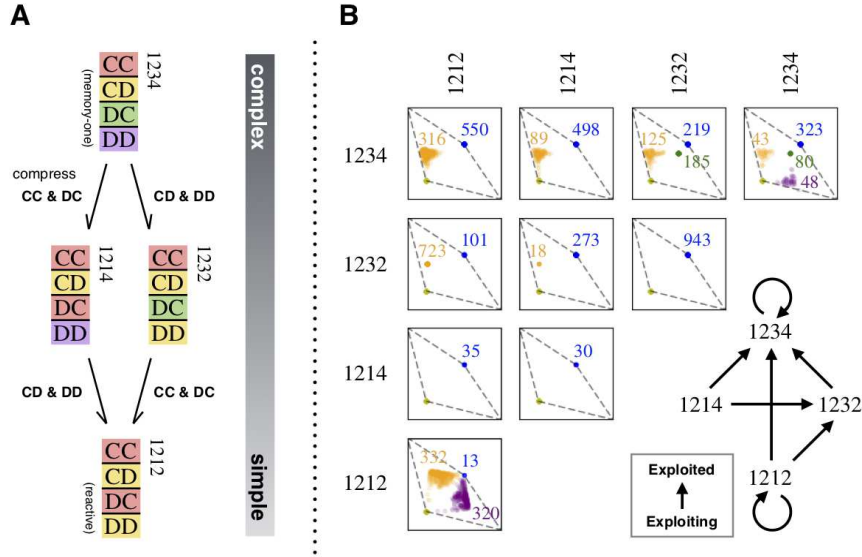}
\caption{{\bf (A)} Extension of classes for the reference to self-action. The 1234 and 1212 classes are equivalent to the previous memory-one and reactive classes, as shown in Fig.~\ref{F01}. By contrast, 1214 and 1232 classes were defined as intermediate classes of 1234 and 1212. The 1214 (1232) class refers to the self-action only when the opponent's previous action is D (C). {\bf (B)} Results of matches among the four classes. The horizontal axis indicates the payoffs of 1234, 1232, 1214, and 1212 from the top to bottom panels, respectively. The vertical axis indicates payoffs of 1212, 1214, 1232, and 1234 in the panels from left to right, respectively. The lower left panel shows the direction of possible exploitation between these classes.}
\label{F06}
\end{center}
\end{figure}

Before analyzing the matches, we labeled these new classes according to Fig.~\ref{F06}. Here, we renamed the memory-one class as the ``1234'' class because the class distinguishes all of CC (1), CD (2), DC (3), and DD (4), and has four strategy variables $x_1$, $x_2$, $x_3$, and $x_4$. Further, the reactive class uses two variables $x_{13}$ and $x_{24}$, as DC (3) and CC (1) represent one class, and DD (4) and CD (2) represent another. Thus, we renamed it as the ``1212'' class. In the same way, the newly defined classes are renamed as ``1214'' and ``1232''; the former (latter) combines CC and DC (CD and DD). Among these four classes, complexity can be introduced as the degree to which the self-action is referred to. Fig.~\ref{F06}-(A) shows the ordering of this complexity. 1234 and 1212 are the most complex and simple of the four classes, respectively, whereas 1214 and 1232 lie between the two.

The outcomes of all possible 10 matches between any pairs of the 1234, 1232, 1214, and 1212 classes are shown in Fig.~\ref{F06}-(B). As shown in these figures, 1212, which is the simplest class, can exploit 1212, 1232, and 1234. Class 1214 exploits 1232 and 1234. Class 1232 can exploit 1234. These exploitative relationships are summarized in the bottom right panel. Interestingly, the results show that the simpler classes generally exploit the more complex classes, but the reverse never occurs. Thus, the reference to self-action generally leads to generosity toward the opponent's defection and accepts exploitation.

Theoretically, there are 15 types of classes for choosing strategies that use the information of the self and the opponent's actions up to the previous turn. The results of the matches over these 15 classes are shown in the Supplementary Material. The ordering of complexity can be defined over these 15 classes. Although, the complexity among these 15 classes does not always indicate the degree of reference to self-action. Interestingly, the one-sided exploitation by a complex class of a simple one is not observed, except for a single case, i.e., the match between the 1232 and 1131 classes.

\subsection{Generality over different payoff matrices}
\label{S5-2}
So far, we have studied the PD game with the standard score $(T,R,P,S)=(5,3,1,0)$. The above results shows that one-sided exploitation by simple classes of complex classes is valid if the payoff satisfies $T-R-P+S>0$. This condition is called the ``submodular'' PD \cite{Dixit2003,Takahashi2010}; the summation of asymmetric payoffs, $T+S$, is larger than the summation of symmetric payoffs, $R+P$.

When the payoff matrix does not satisfy the submodularity, mutual learning does not generate an exploitative relationship. In all the matches among 1234, 1214, 1232, and 1212 classes, the players achieve either mutual cooperation or mutual defection, depending on the initial condition.

\section{Summary and discussion}
\label{S6}
In this study, we investigated how players' payoffs after learning depend on the complexity of their strategies, that is, the degree of reference to previous actions. By extending the coupled replicator model for learning, we formulated an adaptive optimization strategy by learning previous actions. By focusing on a reactive strategy in which the player refers only to the last action of the opponent and a memory-one class in which the player refers to both their own last action and that of the opponent, we uncovered that the latter, which has more information and includes the former, will be exploited by the former, independent of the initial state. Here, both the strategies of the latter (exploited) and former (exploiting) players permanently oscillate as in the prey and predator dynamics, whereas the exploit relationship is maintained. The exploiting (exploited) side uses the narrow-minded (generous) TFT when the previous self-action was defection.

Using this definition, a player using a memory-one class has a larger number of choices of strategy than that of the reactive class. The memory-one class includes extortion strategies \cite{Press2012,Hilbe2013,Stewart2013,Adami2013,Szolnoki2014}, where the player has an advantage over their opponent using a fixed strategy, that is, the player receives a larger payoff than the opponent, independent of the opponent's strategy. Given this, it is quite surprising that the reactive class unilaterally exploits the memory-one class after mutual learning, regardless of the initial strategy. The results show that even if there are possible advantages to choices of strategy, the player may not realize them through learning if the opponent's strategy continues to change. Here, we demonstrated that learning with reference to self-action makes the player generous toward an opponent's defection, as the unknown way to acquire the generosity \cite{Stewart2013,Delton2011,Feng2017,Kurokawa2018}. In this way, learning to obtain a higher payoff with more information counterintuitively results in a poorer payoff than the opponent, who learns with limited information.

It is common for a player to change their next choice depending on past choices. As already seen in reactive and memory-one classes, it is common for a player to change the next choice of behavior depending on their own or opponent's choice. As briefly mentioned in \S~\ref{S5-1}, our formulation can be extended to reference arbitrary information. For instance, we can assume a memory-$n$ strategy, which refers to actions in more than one previous round. Even the two-memory strategy is quite different from the memory-one strategy: The player can use a greater variety of strategies, such as tit-for-2-tat \cite{Boyd1987}. It has been discussed that a reference to multi-memory generates cooperation more efficiently \cite{Hilbe2017}. Our model could be extended to study whether this is true under mutual learning, or whether the player with more information would exploit their opponent or be exploited.

Game theory is often relevant in explaining characteristic human behaviors. The advantage of the TFT strategy indicates poetic justice in human nature. However, humans also reflect on their past behavior. For instance, they could be motivated to perform beneficial actions toward others after betraying others. This study supports how such behavior emerges or is preserved through learning. This, however, provides room for exploitation. Indeed, Eq.~(\ref{E5-016}) shows that the player with reference to their own previous actions becomes generous toward the opponent's defection after the player defected in the previous round. Ironically, this can be exploited.

\section*{Acknowledge}
The authors would like to thank E. Akiyama, and H. Ohtsuki for useful discussions. This research was partially supported by JPSJ KAKENHI Grant Number JP18J13333, and JP21J01393.

\section*{Bibliography}

\newpage

\setcounter{section}{0}
\setcounter{figure}{0}
\renewcommand{\figurename}{FIG. S}
\renewcommand{\thesection}{S\arabic{section}}

\begin{center}
{\LARGE {\bf Supplementary Material}}
\end{center}

\section{General formulation of classes}
In \S~1 of our main manuscript, we defined the previous memory-one class and reactive class. Moreover, in \S~4, we renamed these classes as 1234 and 1212 ones, and additionally defined 1214 and 1232 ones. In general, however, there are 15 classes in total when the player can refer to only the previous one game at the maximum. This section provides the general extension of such classes.

\begin{figure}[H]
\begin{center}
\includegraphics[width=0.6\linewidth]{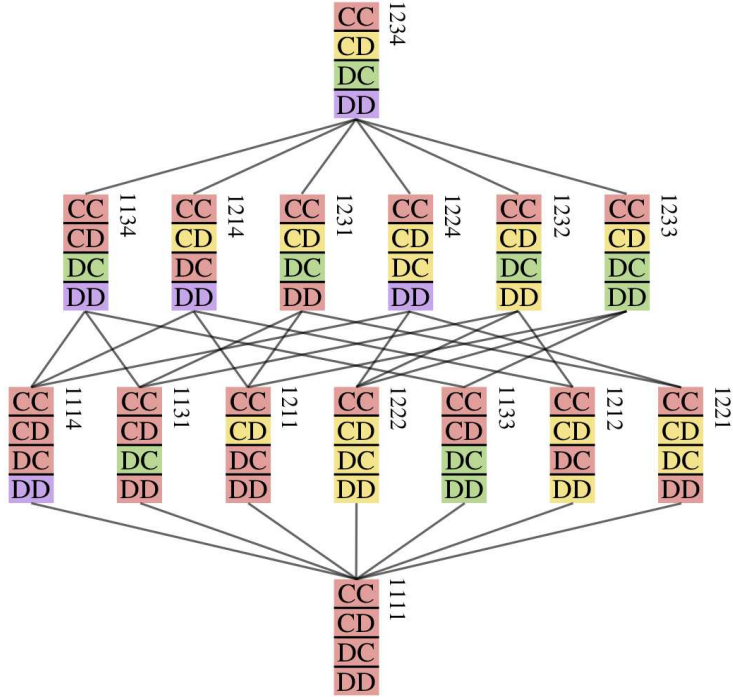}
\caption{Schematic diagram of possible 15 classes and the complexity among them. Each colored bar represents one class of strategies, which consists of four sets of actions, i.e., CC, CD, DC, and DD. When the player reflects the same action to different sets of actions, the sets are painted in the same color. Red, yellow, green and purple indicate numbers 1, 2, 3, and 4. Solid black lines represent the order of complexity between the connected bars. The upper bar is more complex than the connected lower bar.}
\label{FS01}
\end{center}
\end{figure}

Recall that there are four kinds of results in a single game, CC, CD, DC, and DD, as the left and right indices represents the actions of self and other. Then, we assigned numbers to these pairs of actions from 1 to 4, respectively. Then, we assign a number 1234 to the memory-one class, because it distinguishes all of 1 (CC), 2 (CD), 3 (DC), and 4 (DD) and can cooperate in different probabilities depending on the observed action. Next, we consider the reactive class. This reactive class only refers to the other's action, so the player cannot distinguish 3 (DC) with 1 (CC) and 4 (DD) with 2 (CD). When we compress multiple pieces of information, we replace ones with the larger numbers (i.e., 3 and 4) as ones with the smaller numbers (i.e., 1 and 2). Here, 3 and 4 are replaced by 1 and 2 respectively, so that the reactive class is coded as 1212. In the same way, we can define 1134, 1214, 1231, 1224, 1232, 1233, 1114, 1131, 1133, 1211, 1222, and 1221 classes in all. Fig.~S\ref{FS01} shows the schematic diagram of these 15 classes.

In addition, recall the definition of complexity of class. When a class includes all the strategies of another class, the former one is defined to be more complex than the latter one. Fig.~S\ref{FS01} also shows all the relationships of complexity among the classes.

\section{Mutual learning in general classes of strategies}
Before analyzing the games of these 15 classes, we omit 1133 and 1111 classes from the analysis. These two classes do not refer to the other's action. In other words, they do not distinguish ${\rm X_1C}$ with ${\rm X_2D}$ for all ${\rm X_1,X_2}\in\{{\rm C,D}\}^2$. From the rule of prisoner's dilemma, the other player learns to choose D deterministically against 1133 or 1111 classes. Thus, they have no other equilibrium than DD.

\begin{figure}[H]
\begin{center}
\includegraphics[width=0.9\linewidth]{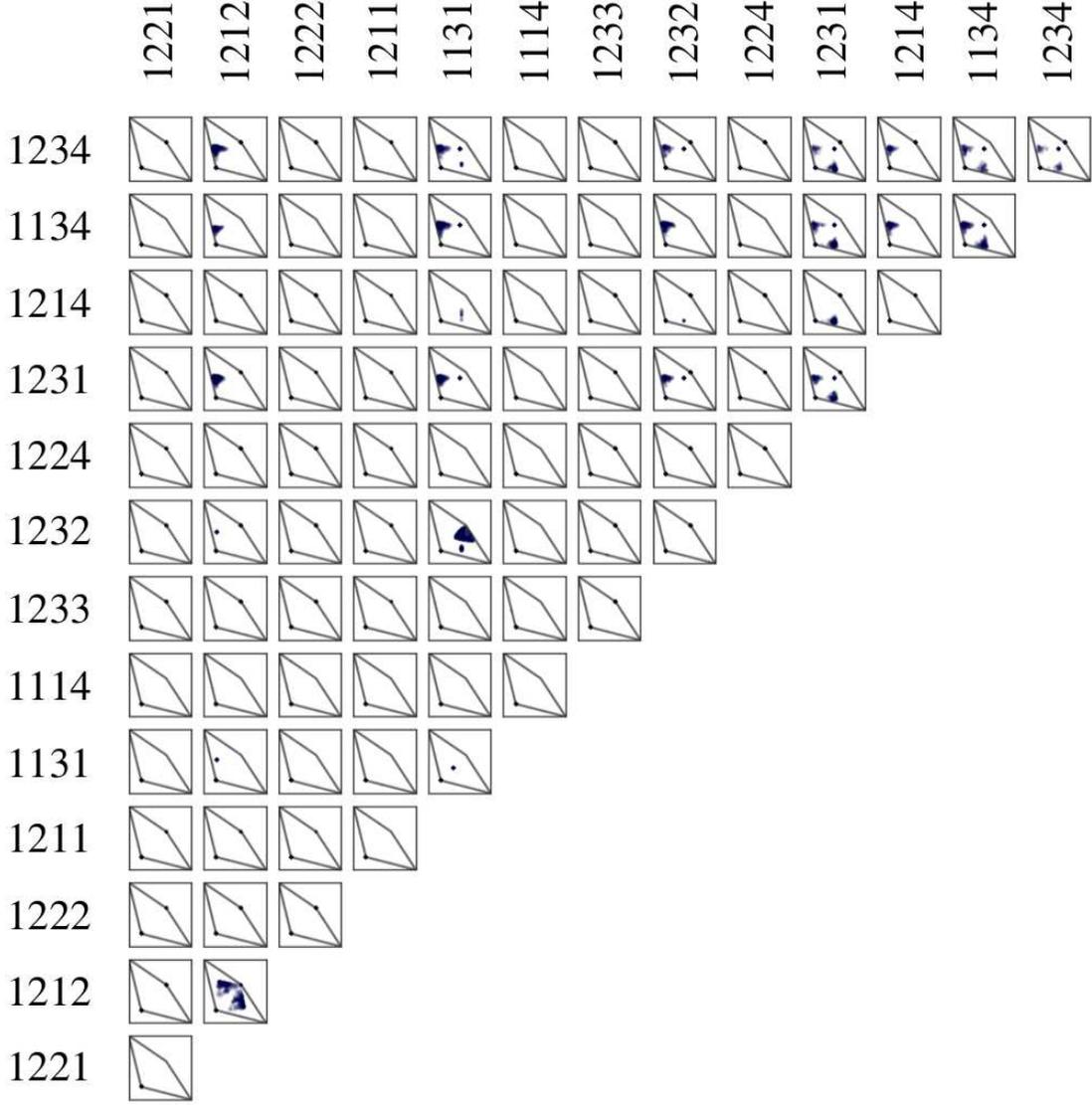}
\caption{The possible equilibrium states in all the pairs of 13 classes of strategies. As the combination of 13 classes, there are $13+13\times 12/2=91$ panels of equilibrium states. In each panel, X (Y) axis indicates the equilibrium payoff of vertical (horizontal) class. Each panel has 10000 samples of equilibrium states in mutual learning with payoff $(T,R,P,S)=(5,3,1,0)$, even though some of equilibria are not plotted.}
\label{FS02}
\end{center}
\end{figure}

Fig.~S\ref{FS02} shows the equilibrium of mutual learning for still existing 13 classes. From this figure, we see a variety of equilibria which include various degrees of cooperation and exploitation. Here, note that there are several kinds of oscillation in the payoffs of both the players, as similar to the case of 1234 v.s. 1234 in the main manuscript. Fig.~S\ref{FS02} shows that the same oscillation is frequently seen in the game of 1234, 1134, 1214, 1231, 1232, 1131, and 1212. Another oscillatory state is seen in the case of 1232 v.s. 1131 on the upper right triangle. All the other states than these two types of oscillation exist as the fixed point in equilibrium of learning.

Then, we give several remarks on computational methods on mutual leanring. First, we gives a constant bound on the stochastic strategies $\epsilon\le x_i\le 1-\epsilon$ with $\epsilon=10^{-4}$ in the computation of learning. This is to avoid false convergences into CC (i.e., $u_{\rm e}=v_{\rm e}=3$) if CC is saddle point. Second, we remove equilibria on $u_{\rm e}=3, 1<v_{\rm e}<3$ or vise versa on several panels. This is because these equilibria only exist on the condition of $R=(T+P)/2$.

\begin{figure}[H]
\begin{center}
\includegraphics[width=0.8\linewidth]{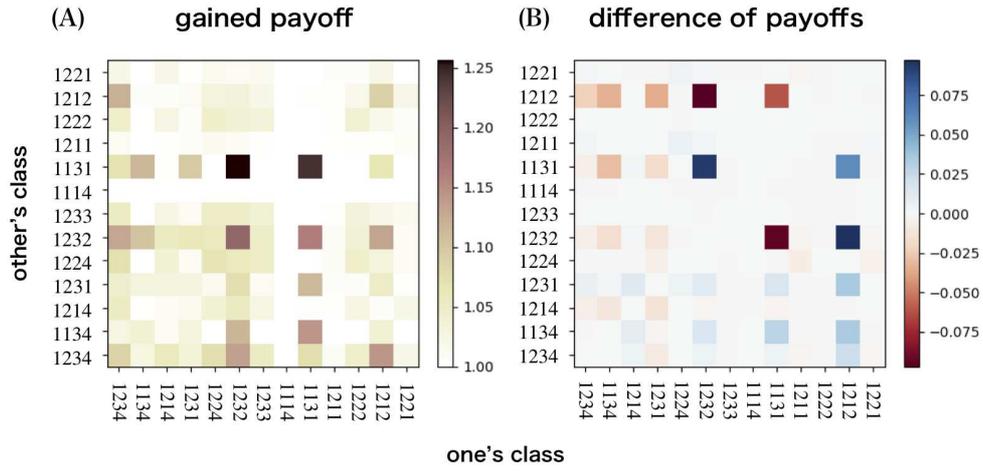}
\caption{Statistical results of Fig.~S\ref{FS02}. Panel (A) shows the payoff of class of horizontal axis number in the ensemble of games against the class of vertical axis. The darker color represents that the class of vertical axis gains a larger payoff. Panel (B) shows the difference in payoff from the horizontal class to the vertical one. The blue (red) color means that the class of horizontal axis exploits (is exploited from) that of vertical axis.}
\label{FS03}
\end{center}
\end{figure}

Fig.~S\ref{FS03} gives a statistical analysis corresponding to Fig.~S\ref{FS02}. The figure (A) shows the payoff of each class obtained statistically from the numerous ensembles. In principle, a learning player receives at least $P(=1)$ amount of payoff in equilibrium. Thus, the difference from such minimal payoff $1$ represents the class's surplus benefit. Interestingly, we see that besides previous 1234 (i.e., memory-one) and 1212 (i.e., reactive) classes, 1232 and 1131 classes achieve high scores in statistics. The figure (B) shows the difference of payoffs between the two players. Interestingly, no exploitation from the simpler class to the more complex one is seen with an exception of 1232 v.s. 1131.


\begin{thebibliography}{99}


\bibitem{vonNeumann2007}
von Neumann J, \& Morgenstern O (2007). {\sl Theory of games and economic behavior (commemorative edition)}. Princeton university press.

\bibitem{Premack1978}
Premack D, \& Woodruff G (1978) Does the chimpanzee have a theory of mind?. {\sl Behavioral and brain sciences}, 1(4), 515-526.

\bibitem{Saxe2007}
Saxe R, Baron-Cohen S, eds (2007) {\sl Theory of Mind: A Special Issue of Social Neuroscience} (Psychology Press, London).

\bibitem{Lurz2011}
Lurz RW (2011) {\sl Mindreading animals: the debate over what animals know about other minds} (MIT Press, Cambridge, MA).

\bibitem{Han2011}
Han TA, Moniz PL, \& Santos FC (2011). Intention recognition promotes the emergence of cooperation. {\sl Adaptive Behavior}, 19(4), 264-279.

\bibitem{Han2015}
Han TA, et al. (2015). Synergy between intention recognition and commitments in cooperation dilemmas. {\sl Scientific reports}, 5(1), 1-7.

\bibitem{Moniz2015}
Moniz PL, Santos FC, \& Lenaerts T (2015). Emergence of cooperation via intention recognition, commitment and apology–a research summary. {\sl Ai Communications}, 28(4), 709-715.

\bibitem{Fujimoto2019a}
Fujimoto Y, \& Kaneko K (2019) Functional dynamic by intention recognition in iterated games. {\sl New Journal of Physics}, 21(2), 023025.

\bibitem{Axelrod1984}
Axelrod R (1984), {\sl The Evolution of Cooperation}, Basic Books, US.

\bibitem{Axelrod1981}
Axelrod R, \& Hamilton WD (1981) The evolution of cooperation. {\sl science}, 211(4489), 1390-1396.

\bibitem{Axelrod1988}
Axelrod R, \& Dion D (1988) The further evolution of cooperation. {\sl Science}, 242(4884), 1385-1390.

\bibitem{Nowak1990a}
Nowak M (1990). Stochastic strategies in the prisoner's dilemma. {\sl Theoretical population biology}, 38(1), 93-112.

\bibitem{Nowak1990b}
Nowak M, \& Sigmund K (1990) The evolution of stochastic strategies in the prisoner's dilemma. {\sl Acta Applicandae Mathematica}, 20(3), 247-265.

\bibitem{Nowak1992}
Nowak MA, \& Sigmund K (1992) Tit for tat in heterogeneous populations. {\sl Nature}, 355(6357), 250-253.

\bibitem{Imhof2010}
Imhof LA, \& Nowak MA (2010) Stochastic evolutionary dynamics of direct reciprocity. {\sl Proceedings of the Royal Society B: Biological Sciences}, 277(1680), 463-468.

\bibitem{Zhang2011}
Zhang J, et al (2011) Resolution of the stochastic strategy spatial prisoner's dilemma by means of particle swarm optimization. {\sl PloS one}, 6(7), e21787.

\bibitem{Baek2016}
Baek SK, et al (2016) Comparing reactive and memory-one strategies of direct reciprocity. {\sl Scientific reports}, 6(1), 1-13.

\bibitem{Nowak1993}
Nowak M, \& Sigmund K (1993) A strategy of win-stay, lose-shift that outperforms tit-for-tat in the Prisoner's Dilemma game. {\sl Nature}, 364(6432), 56-58.

\bibitem{Brauchli1999}
Brauchli K, Killingback T, \& Doebeli M (1999) Evolution of cooperation in spatially structured populations. {\sl Journal of theoretical biology}, 200(4), 405-417.

\bibitem{Kraines2000}
Kraines DP, \& Kraines VY (2000) Natural selection of memory-one strategies for the iterated prisoner's dilemma. {\sl Journal of Theoretical Biology}, 203(4), 335-355.

\bibitem{Iliopoulos2010}
Iliopoulos D, Hintze A, \& Adami C (2010) Critical dynamics in the evolution of stochastic strategies for the iterated prisoner's dilemma. {\sl PLoS computational biology}, 6(10), e1000948.

\bibitem{Stewart2014}
Stewart, AJ, \& Plotkin JB (2014) Collapse of cooperation in evolving games. {\sl Proceedings of the National Academy of Sciences}, 111(49), 17558-17563.

\bibitem{Hilbe2017}
Hilbe C, Martinez-Vaquero LA, Chatterjee K, \& Nowak MA (2017) Memory-n strategies of direct reciprocity. {\sl Proceedings of the National Academy of Sciences}, 114(18), 4715-4720.

\bibitem{Traulsen2006}
Traulsen A, \& Nowak MA (2006) Evolution of cooperation by multilevel selection. {\sl Proceedings of the National Academy of Sciences}, 103(29), 10952-10955.

\bibitem{Wang2013}
Wang Z, Szolnoki A, \& Perc M (2013) Interdependent network reciprocity in evolutionary games. {\sl Scientific reports}, 3(1), 1-7.

\bibitem{Jin2014}
Jin Q, Wang L, Xia CY, \& Wang Z (2014). Spontaneous symmetry breaking in interdependent networked game. {\sl Scientific reports}, 4, 4095.

\bibitem{Xia2018}
Xia C, et al (2018) Doubly effects of information sharing on interdependent network reciprocity. {\sl New Journal of Physics}, 20(7), 075005.

\bibitem{Liu2019}
Liu Y, et al (2019) Swarm intelligence inspired cooperation promotion and symmetry breaking in interdependent networked game. {\sl Chaos: An Interdisciplinary Journal of Nonlinear Science}, 29(4), 043101.

\bibitem{Takesue2021}
Takesue H (2021) Symmetry breaking in the prisoner's dilemma on two-layer dynamic multiplex networks. {\sl Applied Mathematics and Computation}, 388, 125543.

\bibitem{Macy1991}
Macy MW (1991) Learning to cooperate: Stochastic and tacit collusion in social exchange. {\sl American Journal of Sociology}, 97(3), 808-843.

\bibitem{Macy1998}
Macy MW, \& Skvoretz J (1998) The evolution of trust and cooperation between strangers: A computational model. {\sl American Sociological Review}, 638-660.

\bibitem{Macy2002a}
Macy MW, \& Flache A (2002) Learning dynamics in social dilemmas. {\sl Proceedings of the National Academy of Sciences}, 99(suppl 3), 7229-7236.

\bibitem{Macy2002b}
Macy MW, \& Sato Y (2002) Trust, cooperation, and market formation in the US and Japan. {\sl Proceedings of the National Academy of Sciences}, 99(suppl 3), 7214-7220.

\bibitem{Fujimoto2019b}
Fujimoto Y, \& Kaneko K (2019) Emergence of exploitation as symmetry breaking in iterated prisoner's dilemma. {\sl Physical Review Research}, 1(3), 033077.

\bibitem{Sandholm1996}
Sandholm TW, \& Crites RH (1996) Multiagent reinforcement learning in the iterated prisoner's dilemma. {\sl Biosystems}, 37(1-2), 147-166.

\bibitem{Taiji1999}
Taiji M, \& Ikegami T (1999). Dynamics of internal models in game players. {\sl Physica D: Nonlinear Phenomena}, 134(2), 253-266.

\bibitem{Masuda2009}
Masuda N, \& Ohtsuki H (2009) A theoretical analysis of temporal difference learning in the iterated prisoner's dilemma game. {\sl Bulletin of mathematical biology}, 71(8), 1818-1850.

\bibitem{Leibo2017}
Leibo JZ, et al (2017) Multi-agent reinforcement learning in sequential social dilemmas. {\sl arXiv preprint} arXiv:1702.03037.

\bibitem{Sigmund1998}
Hofbauer J, \& Sigmund K (1998). {\sl Evolutionary games and population dynamics}. Cambridge university press.

\bibitem{Borgers1997}
Börgers T, \& Sarin R (1997) Learning through reinforcement and replicator dynamics. {\sl Journal of Economic Theory}, 77(1), 1-14.

\bibitem{Sato2002}
Sato Y, Akiyama E, \& Farmer JD (2002) Chaos in learning a simple two-person game. {\sl Proceedings of the National Academy of Sciences}, 99(7), 4748-4751.

\bibitem{Sato2003}
Sato Y, \& Crutchfield JP (2003) Coupled replicator equations for the dynamics of learning in multiagent systems. {\sl Physical Review E}, 67(1), 015206.

\bibitem{Posch1999}
Posch M (1999) Win–stay, lose–shift strategies for repeated games—memory length, aspiration levels and noise. {\sl Journal of theoretical biology}, 198(2), 183-195.

\bibitem{Imhof2007}
Imhof LA, Fudenberg D, \& Nowak MA (2007) Tit-for-tat or win-stay, lose-shift?. {\sl Journal of theoretical biology}, 247(3), 574-580.

\bibitem{Amaral2016}
Amaral MA, et al (2016) Stochastic win-stay-lose-shift strategy with dynamic aspirations in evolutionary social dilemmas. {\sl Physical Review E}, 94(3), 032317.

\bibitem{Nash1950}
Nash JF (1950) Equilibrium points in n-person games. {\sl Proceedings of the national academy of sciences}, 36(1), 48-49.

\bibitem{Dixit2003}
Dixit A (2003) On modes of economic governance. {\sl Econometrica}, 71(2), 449-481.

\bibitem{Takahashi2010}
Takahashi S (2010) Community enforcement when players observe partners' past play. {\sl Journal of Economic Theory}, 145(1), 42-62.

\bibitem{Press2012}
Press WH, \& Dyson FJ (2012) Iterated Prisoner's Dilemma contains strategies that dominate any evolutionary opponent. {\sl Proceedings of the National Academy of Sciences}, 109(26), 10409-10413.

\bibitem{Hilbe2013}
Hilbe C, Nowak MA, \& Sigmund K (2013) Evolution of extortion in iterated prisoner's dilemma games. {\sl Proceedings of the National Academy of Sciences}, 110(17), 6913-6918.

\bibitem{Stewart2013}
Stewart AJ, \& Plotkin JB (2013) From extortion to generosity, evolution in the iterated prisoner's dilemma. {\sl Proceedings of the National Academy of Sciences}, 110(38), 15348-15353.

\bibitem{Adami2013}
Adami C, \& Hintze A (2013) Evolutionary instability of zero-determinant strategies demonstrates that winning is not everything. {\sl Nature communications}, 4(1), 1-8.

\bibitem{Szolnoki2014}
Szolnoki A, \& Perc M (2014) Evolution of extortion in structured populations. {\sl Physical Review E}, 89(2), 022804.

\bibitem{Delton2011}
Delton AW, et al (2011) Evolution of direct reciprocity under uncertainty can explain human generosity in one-shot encounters. {\sl Proceedings of the National Academy of Sciences}, 108(32), 13335-13340.

\bibitem{Feng2017}
Feng X, Zhang Y, \& Wang L (2017) Evolution of stinginess and generosity in finite populations. {\sl Journal of theoretical biology}, 421, 71-80.

\bibitem{Kurokawa2018}
Kurokawa S, Wakano JY, \& Ihara Y (2018) Evolution of groupwise cooperation: generosity, paradoxical behavior, and non-linear payoff functions. {\sl Games}, 9(4), 100.

\bibitem{Boyd1987}
Boyd R, \& Lorberbaum JP (1987) No pure strategy is evolutionarily stable in the repeated prisoner's dilemma game. {\sl Nature}, 327(6117), 58-59.



\end{thebibliography}
\end{document}